\documentstyle[12pt]{article} 
\def\l{\lambda}
\def\s{\sigma}
\def\ap{a^{\dagger}}
\def\Q{Q^{\dagger}}
\def\H{{\cal H}}
\def\e{{\rm e}}
\def\i{{\rm i}}
\def\be{\begin{equation}}
\def\ee{\end{equation}}
\def\r{\rangle}
\def\Z{\rm Z}
\def\R{\rm R}
\def\C{\rm C}
\begin{document}
\title{$C_{\l}$-extended oscillator algebra and
parasupersymmetric quantum mechanics
\footnote{Presented at the 7th Colloquium ``Quantum Groups and 
Integrable Systems", Prague, 18--20 June 1998}}
\author{C. Quesne \thanks{Directeur de recherches FNRS. E-mail:
cquesne@ulb.ac.be}\ , N. Vansteenkiste \thanks{E:mail:
nvsteen@ulb.ac.be}\\
{\small \sl Physique Nucl\'eaire Th\'eorique et Physique
Math\'ematique,  Universit\'e Libre de Bruxelles,} \\
{\small \sl Campus de la
Plaine CP229, Boulevard~du Triomphe, B-1050 Brussels, Belgium}}
\maketitle

\begin{abstract} The $C_{\l}$-extended oscillator algebra is generated
by $ \{ {\bf 1}, a, \ap, N , T \} $, where $ T $ is the generator of the
cyclic group $C_{\l}$ of order $\l$. It can be realized as a
generalized deformed oscillator algebra (GDOA). Its unirreps can thus be
easily exhibited using the representation theory of GDOAs and their carrier
space show a $\Z_{\l}$-grading structure. Within its infinite-dimensional 
Fock space representation, this algebra provides a bosonization of
parasupersymmetric quantum mechanics of order $p= \l - 1 $. \end{abstract}

\section{Introduction}
Deformations of the oscillator algebra have found a lot of applications to 
physical problems. Some of them are the description of systems with 
non-standard statistics, algebraic realizations of some one-dimensional 
exactly solvable potentials ({\em e.g.} the Calogero model), the use of 
their coherent and squeezed states to describe experimental properties 
of light ({\em e.g.} in a Kerr medium).

Supersymmetric quantum mechanics (SSQM) and the concept of shape invariance 
have proved very useful for generating exactly solvable potentials. SSQM 
establishes also a nice symmetry between bosons and fermions, which has been 
enlarged to symmetry between bosons and parafermions to give rise to 
parasupersymmetric quantum mechanics (PSSQM).

The purpose of this communication is to introduce a new oscillator algebra, 
the $C_{\l}$-extended oscillator algebra, where $C_{\l}$ is a cyclic group of 
order $\l$, and to show a possible bosonization of PSSQM using its elements.
This is one of the prospected results announced in Ref.~\cite{clext}.

\section{Definitions, Fock-space representation and bosonic Hamiltonian}
The Calogero-Vasiliev algebra \cite{vas} is generated by the operators $
\{ {\bf 1}, a , \ap , N , K \}$ satisfying
\begin{eqnarray}
  \left[ a , \ap \right] & = & {\bf 1} + \nu K , \qquad \nu \in \R , 
    \qquad K^2 = {\bf 1} ,   \label{r1va} \\
  \left[ N , a \right] & = & - a , \qquad \left[ N , \ap \right]  =  \ap, 
    \qquad [N, K] = 0 , \label{r2va} \\
  a K & = & - K a , \qquad \ap K = - K \ap ,  \label{r3va} \\
  N^{\dagger} & = & N , \qquad \left(\ap\right)^{\dagger} = a , 
    \qquad K^{\dagger} = K^{-1}  \label{r4va} .
\end{eqnarray}
Klein operator $ K $ can be interpreted as the generator of the symmetric
goup $S_2$. This led to a generalization of the Calogero-Vasiliev algebra
known as the $ S_N $-extended oscillator algebra, used in the study of
the N-particle Calogero problem.

On the other hand, one can see $ K $ as the generator of the cyclic
group $C_2$ of order 2. As proposed in Ref.~\cite{clext}, this leads 
to another generalization, the
$C_{\l}$-extended oscillator algebra ($\l \in \{ 2,3,4,\ldots \}$)
generated by $ \{ {\bf 1}, a , \ap , N , T \}$, satisfying the relations
\begin{eqnarray}
  \left[ a , \ap \right] & = & {\bf 1} + \sum_{\mu=1}^{\l-1} 
                                          \kappa_{\mu} T^{\mu} , 
   \qquad \kappa_{\mu} \in \C , \qquad T^{\l} = {\bf 1} ,    \label{r1cla} \\
  \left[ N , a \right] & = & - a , \qquad \left[ N , \ap \right]  =  \ap, 
   \qquad [N, T] = 0 , \label{r2cla} \\
  a T & = & \e^{\i 2 \pi / \l} \, T a , \qquad 
    \ap T = \e^{- \i 2 \pi / \l} \, T \ap , 
   \label{r3cla} \\
  N^{\dagger} & = & N , \qquad \left(\ap\right)^{\dagger} = a , 
   \qquad T^{\dagger} = T^{-1} .
   \label{r4cla}  
\end{eqnarray}
Here $ T $ is the generator of a unitary representation of the cyclic
group of order $\l$, $C_{\l} = \{ {\bf 1} , T , T^2, \ldots ,T^{\l -1}
\}$, and the $ \l -1$ complex parameters $\kappa_{\mu}$ are restricted
by the condition $\kappa_{\mu}^{\ast} = \kappa_{\l - \mu}$  to preserve 
the Hermiticity
properties of (\ref{r1cla}), so that there are only $\l -1$ independant
real parameters in the definition of the algebra. It is worth noting
that the annihilation and creation operators (resp. $a$ and $\ap$)
quommute with $T$, with $q$ being a $\l$-th root of unity. For $\l=2$, 
this reduces
to the anticommutation (\ref{r3va}) of $a$ and $\ap$ with $K$.

$C_{\l}$ possesses $\l$ inequivalent one-dimensionnal unirreps labelled by
$\mu$, and such that $\Gamma^{\mu}\left(T^\nu\right)=\exp(\i 2 \pi \mu
\nu / \l)$, $\mu , \nu=0,1,2,\ldots,\l-1$. The projector $P_{\mu}$ on
the carrier space of the $\mu$-indexed unirrep takes the form
\be
  P_{\mu} = \frac{1}{\l} \sum_{\nu=0}^{\l-1} \Bigl(\Gamma^{\mu} \left(
  T^{\nu}\right) \Bigr)^{\ast} T^{\nu} = \frac{1}{\l} \sum_{\nu=0}^{\l-1}
  \e^{-\i 2\pi \mu\nu/\l}\, T^{\nu}.   \label{proj}
\ee
It is straightforward to show that these operators satisfy the projection 
operator
defining relations, i.e., $P_{\mu} P_{\nu} = \delta_{\mu,\nu}
P_{\mu}$, and $\sum_{\mu=0}^{\lambda-1} P_{\mu} = {\bf 1}$. Conversely,
$T^{\nu}$ can be expressed in terms of $P_{\mu}$,
\be 
T^{\nu} = \sum_{\mu=0}^{\lambda-1} \e^{\i 2\pi \mu\nu/\lambda} P_{\mu},
\label{invproj}
\ee
so that the $C_{\l}$-extended oscillator algebra can be redefined in
terms of $P_{\mu}$ instead of $T^{\mu}$, 
\begin{eqnarray}
  \left[ a , \ap \right] & = & {\bf 1} + \sum_{\mu=0}^{\l-1} 
                                          \alpha_{\mu} P_{\mu} , 
   \qquad \alpha_{\mu} \in \R , \qquad \sum_{\mu=0}^{\l-1} \alpha_{\mu} = 0 ,  
   \label{r1clap} \\
  \left[ N , a \right] & = & - a , \qquad \left[ N , \ap \right]  =  \ap, 
   \qquad [N, P_{\mu}] = 0 , \label{r2clap} \\
  a P_{\mu} & = & P_{\mu -1} a , \qquad \ap P_{\mu} = P_{\mu +1} \ap , 
   \qquad P_{\mu} P_{\nu} = \delta_{\mu,\nu} P_{\mu} , \label{r3clap} \\
  N^{\dagger} & = & N , \qquad \left(\ap\right)^{\dagger} = a , 
   \qquad \left(P_{\mu}\right)^{\dagger} = P_{\mu} , 
   \qquad \sum_{\mu=0}^{\lambda-1} P_{\mu} = {\bf 1} , \label{r4clap} 
\end{eqnarray}  
where $\mu , \nu = 0,1,2,\ldots,\l-1$, and we use the convention that for any
parameter or operator $X$ involved in a $C_{\l}$-extended oscillator
algebra, $X_y \equiv X_{y \, {\rm mod} \l}$ ({\em  e.g.} 
$P_{-1} = P_{\l-1}$).
The reality and sum conditions on $\alpha_{\mu}$ in (\ref{r1clap}) follow from
their dependence on $\kappa_{\mu}$, i.e., $\alpha_{\mu} =
\sum_{\nu=1}^{\l-1} \exp(\i 2\pi \mu\nu/\l) \kappa_{\nu}$, 
$\mu=0, 1,\ldots,\l-1$. 

$T$ can be realized in terms of the number operator $N$. There exist 
other
realizations (using a function of spin matrices, for instance), 
but here, only that involving $N$ is considered.
Explicitely
\be
  T = \e^{\i 2\pi N/\l}, \qquad P_{\mu} = \frac{1}{\l}
  \sum_{\nu=0}^{\l-1} \e^{\i 2\pi  \nu (N-\mu)/\l}, \qquad \mu = 0, 1,
  \ldots, \l-1 ,     \label{realn}
\ee
assuming the spectrum of $N$ is made of integers only. 
Within this realization, the $C_{\l}$-extended oscillator algebra can be seen 
as
a generalized deformed oscillator algebra 
\linebreak
(GDOA), and its unirreps can be 
developped using the results of Ref.~\cite{rgdoa}.

A GDOA is characterized by an analytic deformation function $G(N)$ such that 
$[a , \ap]_q = G(N)$, and a structure function $F(N)$ vanishing for $N=0$, and
solution of the functional equation $F(N+1) - q F(N) = G(N)$. A GDOA Casimir 
operator is $C= q^{-N} (F(N) - \ap a)$. Here, (\ref{r1clap})
together with (\ref{realn})
give $q=1$, $G(N)= {\bf 1} + \sum_{\mu=0}^{\l-1} \alpha_{\mu} P_{\mu}$, and
\be
F(N)=N+ \sum_{\mu =0}^{\l-1} \beta_{\mu} P_{\mu} , \qquad
 \beta_{\mu} = \sum_{\nu =0}^{\mu-1} \alpha_{\nu} , \qquad \beta_0=0 ,
                             \label{stru}
\ee
with $F(N+1) = N+{\bf1}+ \sum_{\mu =0}^{\l-1} \beta_{\mu+1} P_{\mu}$.

A Fock-space representation is characterized by the existence of a normalized
simultaneous eigenvector $|0\r$ of $N$ and $C$ with both eigenvalues
equal to $0$, and which is destroyed by $a$. In such a unirrep,
$\ap a = F(N)$, $ a \ap =F(N+1)$, and $P_{\mu} |0\r = \delta_{\mu,0} |0\r$.
Its carrier space is spanned by the vectors 
\be
|n\r = {\cal N}_n^{-1/2} \left(\ap\right)^n |0\r, 
                            \qquad n = 0,1, 2, \ldots, d-1, \label{n}
\ee
where $d$ may be finite or infinite, and the normalization coefficients 
are given by 
${\cal N}_n= \prod_{\mu=1}^n F(\mu) = \prod_{\mu=1}^n (\mu + \beta_{\mu})$.
If $|n\r \equiv |k \l + \mu \r$, where $k \in \Z^+$ and 
$\mu = n \, {\rm mod} \l $, does exist,
then the generators of the algebra act on it as
\begin{eqnarray}
N \, |k \l + \mu \r &=& (k \l + \mu ) \, |k \l + \mu \r , \qquad 
  P_{\nu} \, |k \l + \mu \r = \delta_{\nu , \mu} \, |k \l + \mu \r , \\
a \, |k \l + \mu \r &=& \sqrt{k \l + \mu + \beta_{\mu}} \, |k \l + \mu-1 \r ,\\
\ap \, |k \l + \mu \r &=& \sqrt{k \l + \mu +1+ \beta_{\mu+1}} \, 
                                                            |k \l + \mu+1 \r.
  \label{agn}
\end{eqnarray}
The existence condition for $|n\r$ is the unitarity condition of the 
representation,
which is fullfilled if and only if $F(\mu) \geq 0$, 
$\forall \mu = 1,2,\ldots,d$.
There
are thus two types of Fock space representations:
\begin{enumerate}
 \item finite-dimensional unirreps of dimension $d < \l$\\
       $\Leftrightarrow F(\mu) > 0$, $\forall \mu= 1,2,\ldots,d-1$, and 
       $F(d)=0$,
 \item infinite-dimensional bounded from below (BFB) unirrep\\
       $\Leftrightarrow F(\mu) > 0$, or 
       $\beta_{\mu} = \sum_{\nu=0}^{\mu-1} \alpha_{\nu} > - \mu$,
       $\forall \mu= 1,2,\ldots,\l-1$.
\end{enumerate}
From now on, only the bosonic Fock space representation, corresponding to 
the second type, will be considered. 

It is worth noting the Fock-space $\Z_{\l}$-grading: for each 
$\mu = 0,1,\ldots, \l-1 $, let
${\cal F}_{\mu} = \{ |k \l + \mu \r $, $k \in \Z^+ \}$; 
then the entire carrier Fock space is
${\cal F}= \sum_{\mu=0}^{\l-1} \oplus {\cal F}_{\mu}$, and the action 
of the algebra generators on a vector of ${\cal F}$ depends on which 
${\cal F}_{\mu}$ subspace it belongs.

The bosonic oscillator Hamiltonian is defined by
$ H_0 = \frac{1}{2} \{ \ap , a \} = \frac{1}{2} (F(N)+F(N+1))$. By
using (\ref{stru}), $H_0$ can be rewritten as
\be
H_0 = N + \frac{1}{2} {\bf 1} + \sum_{\mu=0}^{\l-1} \gamma_{\mu} P_{\mu},
 \qquad \gamma_{\mu}= \sum_{\nu=0}^{\mu-1} \alpha_{\nu} + 
                                       \alpha_{\mu} /2 
 \qquad (\gamma_0 = \alpha_0 /2 ) , \label{habo}
\ee
and it is obvious that
\be
H_0 \, |n\r = E_{n} |n\r , \qquad {\rm where} \qquad
 E_n = n+\frac{1}{2} + \gamma_{n \, {\rm mod} \l}.  \label{en}
\ee
The spectrum of the bosonic Hamiltonian possesses $\l$ families of 
equally spaced eigenvalues (within each ${\cal F}_{\mu}$, $H_0$ is
harmonic). Because of this very interresting property, the 
$C_{\l}$-extended oscillator algebra provides an algebraic realization
for the recently introduced cyclic shape invariant potentials 
\cite{sukha}. This will not be developed here, for more details see
Refs.~\cite{clext,c3ext}.

\section{Bosonization of (para)supersymmetric quantum \,
mechanics}
Supersymmetric quantum mechanics (SSQM) is characterized by supercharge
operators $Q$, $\Q$ such that
\be 
Q^2=0  , \qquad \{ \Q , Q \} = \H ,  \qquad 
 \left( \Q \right)^{\dagger} = Q  ,
\ee
so that the supersymmetric Hamiltonian $\H$ commutes with the supercharges
\linebreak
($ [ \H , Q ] = 0 $). There exists a realization of this algebra wherein
the supercharges are the product of mutually commuting fermionic (Pauli 
matrices) and bosonic operators. An alternative realization 
of SSQM, without fermionic matrices, has been provided using the 
Calogero-Vasiliev (or $C_2$-extended oscillator) algebra 
generators~\cite{bossqm}. The choice 
$ Q = \ap P_1 $ (so that $\H = \ap a \, P_0 + a \ap \, P_1 $) corresponds 
to an unbroken SSQM (except for the ground state, all the states are two-fold 
degenerate), while $ Q = \ap P_0 $ (and $ \H = a \ap \, P_0 + \ap a \, P_1$) 
describes a broken SSQM (all the states are two-fold degenerate).

SSQM has been generalized to parasupersymmetric quantum mechanics 
(PSSQM) of order $p$ by Khare \cite{khare}:
\begin{eqnarray}
Q^{p+1}& = &0 \, , \qquad Q^n \neq 0 \, , \; \; n=1,2,\cdots ,p  , 
          \label{pssqm1}  \\ 
\left[ \H , Q \right] & = & 0  ,  \label{pssqm2}
\end{eqnarray}
and
\be 
Q^p \Q + Q^{p-1} \Q Q + \cdots + \Q Q^{p} = 2 p \, Q^{p-1} \H \, ,
     \label{pssqm3}
\ee
so that SSQM is PSSQM of order 1. PSSQM can be realized in terms of 
mutually commuting 
parafermionic (matrices) and bosonic operators. A property of PSSQM 
of order $p$ is that the energy levels above the $(p-1)$th one are 
$(p+1)$-fold degenerate. This fact, and the bosonization of SSQM by the 
$C_2$-extended oscillator algebra hint at a possibility of describing PSSQM 
in terms of the generators of the $C_{p+1}$-extended oscillator algebra, 
whose bosonic Hamiltonian possesses $p+1$ families of equally spaced 
eigenvalues, which may be shifted to coincide from some excited state onwards. 

Taking $Q_{\s}= \sum_{\nu=0}^{p} \s_{\nu} \, \ap P_{\nu}$, $ \s_{\nu} \in \C$,
where all the operators belong to the $C_{p+1}$-extended oscillator algebra, 
Eq.(\ref{pssqm1}) is obtained iff for only one $\nu$, $\s_{\nu}=0$. Hence 
good candidates for parasupercharges are
\be
Q_{\eta,\mu} = \sum_{\nu=1}^{p} \eta_{\mu+\nu} \, \ap P_{\mu+\nu} , 
                          \label{pscha}
\ee
with $\mu = 0,1, \ldots , p$, and $ \eta_{\mu+\nu} \in \C_0$.

For a given $\mu$, the ansatz 
\be
\H_{\mu} = \frac{1}{2} \left\{ \ap , a \right\} + \frac{1}{2} 
 \sum_{\nu=0}^{p} r_{\nu,\mu} P_{\nu} , \label{psham}
\ee
where $r_{\nu,\mu} \in \R$, commutes with $Q_{\eta,\mu}$ iff 
$[ \H_{\mu} , \ap P_{\mu+\nu} ] = 0$, $\forall \nu = 1,2, \ldots , p$,
which is equivalent to
\be
\begin{array}{ll}
r_{\mu,\mu} \in  \R ,\\
r_{\mu+\nu,\mu} =  2 + \alpha_{\mu+\nu} + \alpha_{\mu+\nu+1} 
                        + r_{\mu+\nu+1,\mu} , \qquad \nu=1,2,\ldots ,p .
                    \label{pshamcond}
\end{array}
\ee

It now remains to impose the nonlinear relation (\ref{pssqm3}) between 
$Q_{\eta,\mu}$ and $\H_{\mu}$, defined in (\ref{pscha}) and (\ref{psham}),
respectively. Such a relation is valid on every state of the BFB Fock-space
representation of the associated $C_{p+1}$-extended oscillator algebra iff  
\begin{eqnarray}
& \sum_{\nu=0}^{p-1} &  | \eta_{\mu+\nu+1} |^2 = 2p , \label{psnl1} \\
& \sum_{\nu=1}^{p-1} &  | \eta_{\mu+\nu+1} |^2 
 \left( \nu + \sum_{\rho=0}^{\nu-1} \alpha_{\mu+\rho+2} \right)= 
 p \left( 1 +\alpha_{\mu+2} + r_{\mu+2,\mu} \right) .
                    \label{psnl2}
\end{eqnarray}

It can be shown that Khare PSSQM of order $p$ can be bosonized with any
$C_{p+1}$-extended oscillator algebra admitting a BFB Fock representation,
or in other words that Eqs.~(\ref{pshamcond}), (\ref{psnl1}), and (\ref{psnl2})
admit solutions for $\eta_{\mu+\nu}$ and $r_{\nu,\mu}$ for any allowed values
of the algebra parameters. For instance, for the choice 
$ |\eta_{\mu+\nu}|^2 = 2$, $\nu = 1,2, \ldots ,p$, satisfying (\ref{psnl1}),
one finds
\be
r_{\mu+2,\mu} = \frac{1}{p} \left\{ (p-2) \, \alpha_{\mu+2} +
  2 \sum_{\nu=\mu+3}^{\mu+p} (p+\mu-\nu+1) \, \alpha_{\nu} + p(p-2) \right\} , 
                                \label{pssol}
\ee
with the remaining $r_{\nu,\mu}$'s taking definite values in terms of
$r_{\mu+2,\mu}$ and $\alpha_{\rho}$. For $\mu =0$, PSSQM is unbroken,
otherwise it is broken with a $(\mu+1)$-fold degenerate ground state.
All the excited states are $(p+1)$-fold degenerate. For 
$\mu = 0,1,\ldots , p-2$, the ground state energy may be positive, null,
or negative depending on the parameters, whereas for $\mu = p-1$ or $p$,
it is always positive.

In Ref.~\cite{clext}, it was shown that for $p=2$, Beckers-Debergh PSSQM
\cite{bd} can only be bosonized with those $C_3$-extended oscillator algebras
for which $\alpha_{\mu+2}= -1$. In such a case, Khare (or Rubakov-Spiridonov)
PSSQM is simultaneously realized.

\section{Conclusion}
A new deformation of the oscillator algebra, called $C_{\l}$-extended 
oscillator algebra, has been introduced, and seen to be a generalization 
of the Calogero-Vasiliev oscillator algebra. 

It provides an algebraic realization of supersymmetric cyclic shape invariant
potentials. By using its generators, PSSQM of order $\l -1$ can be realized 
without parafermionic matrices.

Some open and very interesting questions concern the possibilities of
representing the algebra in terms of differential operators and providing
it with a Hopf structure. Other investigations related to realizations of 
pseudosupersymmetric and orthosupersymmetric QM are investigated.

\end{document}